\documentclass[11pt,a4paper]{article}
\usepackage{amsmath}
\usepackage{amsfonts}
\usepackage{amssymb}
\usepackage{amsthm}
\usepackage{authblk}
\usepackage{mathrsfs}
\usepackage{dsfont}
\usepackage{paralist}
\usepackage{natbib}
\usepackage{bm}
\usepackage{bbm}
\usepackage{color}
\usepackage{graphicx}
\usepackage{tabularx}
\usepackage[left=3.25cm,right=3.25cm,top=2cm,bottom=2cm,includeheadfoot]{geometry}
\usepackage{comment}
\usepackage{tikz}
\usetikzlibrary{decorations.pathreplacing}
\usetikzlibrary{decorations.pathmorphing}
\usepackage{resizegather}
\DeclareFontFamily{OT1}{pzc}{}
\DeclareFontShape{OT1}{pzc}{m}{it}{<-> s * [1.10] pzcmi7t}{}
\DeclareMathAlphabet{\mathpzc}{OT1}{pzc}{m}{it}
\usepackage[colorlinks=true,linkcolor=blue,citecolor=blue,pdfborder={0 0 0}]{hyperref}
\newtheoremstyle{normal}% name
{2ex}               % Space above, empty = `usual value'
{3ex}               % Space below
{}                  % Body font
{}                  % Indent amount (empty = no indent, \parindent = para indent)
{\bfseries} % Thm head font
{}                  % Punctuation after thm head
{2pt}   % Space after thm head.
{\thmname{#1}\thmnumber{ #2.} \thmnote{(#3)}}% Thm head spec

\newtheoremstyle{italic}% name
{2ex}%      Space above, empty = `usual value'
{3ex}%      Space below
{\itshape}% Body font
{}%         Indent amount (empty = no indent, \parindent = para indent)
{\bfseries} % Thm head font
{}%        Punctuation after thm head
{2pt}%     Space after thm head.
{\thmname{#1}\thmnumber{ #2.} \thmnote{(#3)}}% Thm head spec

\theoremstyle{normal}
\newtheorem{definition}{Definition}[section]
\newtheorem{remark}[definition]{Remark}

\theoremstyle{italic}
\newtheorem{theorem}[definition]{Theorem}
\newtheorem{lemma}[definition]{Lemma}

\newtheorem{corollary}[definition]{Corollary}

\renewcommand{\P}{\mathbb{P}}

\newcommand{\E}{\mathbb{E}}

\DeclareMathAccent{\verywidehat}{\mathord}{largesymbols}{'144}

\newcommand{\D}{\mathbb{D}}

\newcommand{\N}{\mathbb{N}}
\newcommand{\R}{\mathbb{R}}

\newcommand{\al}{\alpha}
\newcommand{\be}{\beta}

\newcommand{\si}{\sigma}

\newcommand{\eps}{\varepsilon}

\newcommand{\om}{\omega}

\newcommand{\vp}{\varphi}

\newcommand{\nN}{\mathcal N}

\newcommand{\lL}{\mathcal L}

\newcommand{\tol}{\stackrel{\lL}{\longrightarrow}}

\newcommand{\pn}{\stackrel{\P}{\longrightarrow}}

\newcommand{\mybar}[1]{\makebox[0pt]{$\phantom{#1}\overline{\phantom{#1}}$}#1}

\DeclareMathOperator{\Var}{Var}

\allowdisplaybreaks[1]

\begin{document}

\title{Testing for unspecified periodicities in binary time series}
\author{Finn Schmidtke}
\author{Mathias Vetter\thanks{corresponding author; postal address: Christian-Albrechts-Universit\"at zu Kiel, Mathematisches Seminar, Heinrich-Hecht-Platz\ 6, 24118 Kiel, Germany; e-mail address: vetter@math.uni-kiel.de.}}
\affil{Christian-Albrechts-Universit\"at zu Kiel}

%\emailAdd{vetter@math.uni-kiel.de}

\maketitle

\begin{abstract}
Given independent random variables $Y_1, \ldots, Y_n$ with $Y_i \in \{0,1\}$ we test the hypothesis whether the underlying success probabilities $p_i$ are constant or whether they are periodic with an unspecified period length of $r \ge 2$. The test relies on an auxiliary integer $d$ which can be chosen arbitrarily, using which a new time series of length $d$ is constructed. For this new time series, the test statistic is derived according to the classical $g$ test by Fisher. Under the null hypothesis of a constant success probability it is shown that the test keeps the level asymptotically, while it has power for most alternatives, i.e.\ typically in the case of $r \ge 3$ and where $r$ and $d$ have common divisors.  
\end{abstract}

\medskip


\smallskip


%\tableofcontents

\section{Introduction}
\def\theequation{1.\arabic{equation}}
\setcounter{equation}{0}

Testing for periodicities in time series is among the most classical problems in statistics. For Gaussian processes, an excellent review can be found in \S 10 of \cite{brodav2006} in which the authors present several tests which can ultimately be traced back to the foundations of the field. Of special interest for us is a test by \cite{fisher1929} in which the null hypothesis of Gaussian white noise is tested against the presence of an additional deterministic component of an unknown frequency. Given i.i.d.\ observations $X_1, \ldots, X_n$ with $X_i \sim \nN(0,\si^2)$ and $q = \lfloor (n-1)/2 \rfloor$, the idea behind Fisher's test (sometimes referred to as Fisher's $g$) is to reject the null hypothesis whenever at least one of the periodograms is large relative to the others. Precisely, the test statistic is given by 
\[
 \xi_q = \frac{\max_{j=1, \ldots, q} I_X(\om_j)}{\sum_{m=1}^q I_X(\om_m)}
\]
where $\om_j = \frac{2 \pi j}n$ denotes the $j$th Fourier frequency and 
\begin{align*}
I_X(\om_j) = \frac 1n \left| \sum_{\ell=1}^n X_\ell \exp(-i \ell \om_j) \right |^2 
\end{align*}
is the associated periodogram. Following the discussion leading to Corollary 10.2.2 in \cite{brodav2006}, critical values can be obtained in accordance with Theorem I.9.2 in \cite{feller1971} via the formula
\begin{align} \label{limxi}
\P(\xi_q^n \ge x) = 1- \sum_{j=0}^q (-1)^j \binom{q}{j} \left(1-jx\right)_+^{q-1} = \sum_{j=1}^q (-1)^{j+1} \binom{q}{j} \left(1-jx\right)_+^{q-1}
\end{align}
where $x_+ = \max\{x,0\}$ stands for the positive part of $x$. Note that Fisher's test is exact, i.e.\ no asymptotics in the sense of $n \longrightarrow \infty$ are assumed. 

Since its invention, Fisher's test has been generalised and improved in several aspects. Typically it is assumed that 
\begin{align} \label{lingen}
X_\ell = \mu + \al \cos(\om \ell) + \be \sin(\om  \ell) + \eps_\ell, \quad \ell=1, \ldots, n,
\end{align}
where $\mu, \al, \be$ and $\om \in (0,\pi)$ are unknown parameters and the $\eps_i$ are i.i.d.\ distributed with mean zero and an unknown variance $\si^2$. The null hypothesis then can be stated as $H_0: \al = \be = 0$, while the alternative becomes $H_1: \al > 0 \text{ or } \be > 0$. To mention a few classical works, we refer to Section 4.3 of \cite{anderson1971} for optimality results of Fisher's test, to \cite{walker1965} and \cite{davmik1999} for extensions to other distributions than normal ones, and to \cite{whittle1952}, \cite{siegel1980} and \cite{quinn1986} for more general alternatives than just a single frequency. 

In this paper, however, we are interested in a binary time series. Our standing assumption is that we observe a sequence $Y_1, Y_2, \ldots$ of independent random variables with values in $\{0,1\}$, where for each $Y_\ell$ it is supposed that there exists a deterministic number $p_\ell = \P(Y_\ell = 1) \in (0,1)$. We assume throughout that the sequence $p_1, p_2, \ldots$ is periodic with period length $r$, i.e.\ $r$ is the minimal integer such that $p_\ell = p_{mr + \ell}$ for all $\ell, m \in \N$. A natural translation of the model behind Fisher's test to the binary world is then to test 
\begin{align} \label{h0h1}
H_0: r = 1 \text{ (constant success probability) against } H_1: r \ge 2 \text{ (true periodicity)}
\end{align}
based on a total of $n$ observations. Note that this model cannot be represented in the general form \eqref{lingen} as the error term needs to be of a very specific form and with a time-varying variance in order to make $Y_\ell$ binary.

Binary time series are a standard way to model the evolution of categorical data over time, and they are used in various fields such as neuroscience (\citealp{renetal2010}), economy (\citealp{belfer2010}), ecology (\citealp{adapir2019}) or geophysics (\citealp{tueetal2017}; Chapter 2 of \citealp{kedfok2002}). Albeit most models are more sophisticated and work with success probabilities which are stochastic processes themselves, even purely deterministic models are used in the literature, for example in volcanology (\citealp{ammnav2003}).

Surprisingly little attention has been devoted to the testing problem (\ref{h0h1}). In his unpublished work, \cite{coakley1991} presents a test of $H_0$ against the specific alternative that the period length $r$ equals a given integer $k$, while \cite{wilnav1997} discuss a Bayesian version of his test. Unsurprisingly, the basic idea behind Coakley's test is that each 
\[
X_j^n= \frac kn \sum_{\ell=1}^{\frac nk} Y_{(\ell-1)k+j}, \quad j=1, \ldots, k,
\]
estimates the unknown probability $p_j$, where it is assumed that $n$ is a multiple of $k$. Under the null hypothesis each $X_j^n$ should be close to the sample mean $\mybar{Y}_n$, and hence his test statistic is based on $$\sum_{j=1}^k (X_j^n - \mybar{Y}_n)^2$$ for which a rescaled version is shown to be asymptotically $\chi^2_{k-1}$-distributed under $H_0$ and asymptotically consistent under the alternative of a period $k$.

Our solution to the testing problem in (\ref{h0h1}) can be regarded as a combination of Coakley's procedure and Fisher's test, and it relies on the choice of an auxiliary parameter $d \in \N$ which serves as the length of an artificial time series to which the original one is transformed in a similar way as in Coakley's work. Ideally, $d$ should be chosen as a number which divides a lot of integers. To this new time series of length $d$, and with the associated $q$, we then apply Fisher's test statistic. Under the null hypothesis we show with the aid of the central limit theorem that asymptotically (as $n \longrightarrow \infty$) it obeys the same law as in (\ref{limxi}). Under the alternative, the test has a large power in a lot of cases, namely in general if $d$ and $r$ have common divisors and $r \ge 3$ .

The remainder of the paper is organised as follows: In Section \ref{sec:main} we present the test statistic and state the main results of our work, which finite sample properties we analyse in Section \ref{sec:sim}. All proofs are gathered in Section \ref{sec:proof}.

\section{Main results} \label{sec:main}
\def\theequation{2.\arabic{equation}}
\setcounter{equation}{0}

As noted above, $Y_1, Y_2, \ldots$ denotes a sequence of independent random variables with values in $\{0,1\}$, where for each $Y_\ell$ it is supposed that there exists a deterministic number $p_\ell = \P(Y_\ell = 1) \in (0,1)$. The sequence $p_1, p_2, \ldots$ is assumed to be periodic with period length $r$. We then choose an auxiliary integer $d\in\mathbb{N}$ which serves as the length of a new time series constructed from the observations $Y_1, Y_2, \ldots, Y_n$. In principle, $d$ can be chosen arbitrarily but for reasons becoming apparent later it is beneficial to have a lot of distinct divisors within $d$. For any $1\leq i\leq d$ we then set
\[
Z_{i}^{n}=\frac{1}{\left\lfloor \frac{n}{d}\right\rfloor}\sum_{k=0}^{\left\lfloor \frac{n}{d}\right\rfloor -1} Y_{i+kd}
\]
where we use the notation $\lfloor x \rfloor$ to denote the largest integer $m$ with $m \le x$. 

The idea behind the proposed test statistic now is as follows: Under the null hypothesis of a constant success probability it is clear that the $Z_{i}^{n}$ are i.i.d., and after a suitable rescaling they are asymptotically normal as well. Hence, for a large $n$ the family $Z_{1}^{n}, \ldots, Z_{d}^{n}$ can be regarded as a time series of length $d$ where the elements are (after an affine transformation) close to i.i.d.\ normal distributions. Recalling the model in (\ref{lingen}), where the null hypothesis corresponds to $\al = \be = 0$, and using that Fisher's test statistic is invariant under certain affine linear transformations, it is plausible that a version of Fisher's test yields critical values which at least asymptotically keep the level $\al$. 

Under the alternative, however, the behaviour of the time series is less easy to figure out. We will see in Theorem \ref{thmei} and Lemma \ref{lemei} below that, if $r \le d$ holds, each $Z_{i}^{n}$ converges in probability to some $e_i$ and that the sequence $e_1, \ldots, e_d$ is itself $r$-periodic. At first sight it hence seems reasonable to regard $Z_{1}^{n}, \ldots, Z_{d}^{n}$ as being close to a deterministic periodic function so that Fisher's test should have power as well. The problem, however, is that $e_1, \ldots, e_d$ is often not just $r$-periodic but actually constant as well. In this case, Fisher's test obviously has no power. This happens in particular when $r$ and $d$ do not have common divisors, but it is not limited to this situation. 

\begin{theorem} \label{thmei}
Let $1 \le i \le d$ be fixed and set
\[
e_{i}=\frac{1}{r} \sum_{k=0}^{r-1} p_{i+kd}
\]
where $r \in \N$ denotes the unknown period of the $p_i$. If $r \le d$, then 
\[
Z_i^n \pn e_i. 
\]
\end{theorem}

\begin{lemma} \label{lemei}
Let $e_{i}$ be as in Theorem \ref{thmei}, $1 \le i \le d$, and let $r \le d$. Then the following claims hold: 
\begin{itemize}
	\item[(a)] The sequence $e_1, \ldots, e_d$ is $r$-periodic. 
%\item [(b)] For any $m \in \N$ we have $e_{i}=\frac{1}{mr}\sum_{k=0}^{mr-1} p_{i+kd}$,  $1 \le i \le d$.
\item [(b)] If $r$ and $d$ do not have common divisors, then the $e_{i}$ are constant
in $i$.
\item [(c)] If $r$ and $d$ do have common divisors, $r<d$, then the $e_{i}$ are in general not constant. 
\end{itemize}
\end{lemma}

\begin{remark} \label{remei}
A key step in the proof of Theorem \ref{thmei} is to show 
\[
\underset{n\rightarrow\infty}{\lim}\E[Z_i^n] = \underset{n\rightarrow\infty}{\lim}\frac{1}{\left\lfloor \frac{n}{d}\right\rfloor }\sum_{k=0}^{\left\lfloor \frac{n}{d}\right\rfloor -1}\mathbb{E}\left[Y_{i+kd}\right]=\frac{1}{r}\sum_{j=0}^{r-1}\E\left[Y_{i+jd}\right] = e_{i}
\]
for any fixed $1 \le i \le d$ and with $r \le d$. In fact, the same proof gives 
\[
\underset{n\rightarrow\infty}{\lim}\Var\left(\sqrt{\left\lfloor \frac{n}{d}\right\rfloor} Z_i^n\right) = \underset{n\rightarrow\infty}{\lim} \frac{1}{\left\lfloor \frac{n}{d}\right\rfloor }\sum_{k=0}^{\left\lfloor \frac{n}{d}\right\rfloor -1}\mathrm{Var}\left(Y_{i+kd}\right) = \frac{1}{r}\sum_{j=0}^{r-1}\mathrm{Var}\left(Y_{i+jd}\right) =: v_i >0
\]
as well which will be helpful later in the proof of Theorem \ref{thmgen}. \qed
\end{remark}

In the following, we are interested in the asymptotic behaviour of the test statistic 
\[
g(Z_1^n, \ldots, Z_d^n) = \frac{\max_{j=1, \ldots, q} I_Z^n(\om_j)}{\sum_{m=1}^q I_Z^n(\om_m)}
\]
where 
\[
I_Z^n(\om_j) = \frac 1d \left| \sum_{\ell=1}^d Z^n_\ell \exp(-i \ell \om_j) \right |^2 
\]
denotes the $j$th periodogram associated with $Z_{1}^{n}, \ldots, Z_{d}^{n}$ and where we have set $q = \lfloor (d-1)/2 \rfloor$, similarly to the introduction. Note that $g$ is not well-defined in general as the denominator vanishes whenever $Z^n = (Z_1^n, \ldots, Z_d^n) \in A$ where the (closed) set $A$ is defined via
\begin{align*}
A=&\left\{\left(x_{1},\ldots,x_{d}\right)\in\mathbb{R}^{d} ~\middle|~ I_{x}(\omega_{1})=\ldots=I_{x}(\omega_{q})=0\right\} \\=& \left\{\left(x_{1},\ldots,x_{d}\right)\in\mathbb{R}^{d} ~\middle|~ \sum_{\ell=1}^d x_\ell \exp(-i \ell \om_j)=0 \text{ for all } 1 \le j \le q\right\}.
\end{align*}
For this reason we will formally set the test statistic to be
\[
f(Z_1^n, \ldots, Z_d^n) = g(Z_1^n, \ldots, Z_d^n) \mathds{1}_{\left\{ \left(Z_{1}^{n},\ldots,Z_{d}^{n}\right)\notin A\right\}}.
\]
It turns out that asymptotic behaviour of this statistic crucially depends on whether $\left(e_{1},\ldots,e_{d}\right)$ is an element of $A$ or not. We start with an important lemma paving the way for the asymptotics of $f(Z_1^n, \ldots, Z_d^n)$ when the limits $\left(e_{1},\ldots,e_{d}\right)$ are in $A$. 

\begin{lemma} \label{lemeiA}
Let $\left(e_{1},\ldots,e_{d}\right) \in A$ %and $\left(Z_{1}^n,\ldots, Z_{d}^n\right) \notin A$
 and $c_n > 0$ be arbitrary. Then 
\[
f\left(Z_{1}^{n},\ldots,Z_{d}^{n}\right) = f\left(c_n(Z_{1}^{n}-e_{1}), \ldots, c_n(Z_{d}^{n}-e_{d})\right).
\]
\end{lemma}

On one hand, Lemma \ref{lemeiA} can be used to derive the asymptotics of $f\left(Z_{1}^{n},\ldots,Z_{d}^{n}\right)$ under the null hypothesis because then $\left(e_{1},\ldots,e_{d}\right) \in A$ holds (see the proof of Corollary \ref{thmnull}), but it also indicates that some alternative hypotheses cannot be detected, in particular if they share the same limiting distribution as under the null or at least a closely related one. Key to the proof of the asymptotics for $f\left(Z_{1}^{n},\ldots,Z_{d}^{n}\right)$ is that each component $c_n(Z_{j}^{n}-e_{j})$ converges in law to a normal distribution if we choose 
\[
c_n = \sqrt{\left\lfloor {n}/{d}\right\rfloor}, 
\]
together with a version of the extended continuous mapping theorem.

\begin{theorem} \label{thmgen}
Let $\left(e_{1},\ldots,e_{d}\right) \in A$ and let $v_i$ be as in Remark \ref{remei}, $1 \le i \le d$. Then 
\[
f\left(Z_{1}^{n},\ldots,Z_{d}^{n}\right) \tol g\left(v_1 N_1,\ldots, v_d N_d\right) \stackrel{d}{=} g\left(N_1, \frac{v_2}{v_1} N_2, \ldots, \frac{v_d}{v_1} N_d\right)
\]
where $N_1, \ldots, N_d$ are i.i.d.\ standard normal distributed random variables. 
\end{theorem}

\begin{corollary} \label{thmnull}
Suppose that the null hypothesis of a constant success probability holds. 
\begin{itemize}
	\item[(a)] Then we have 
\[
f\left(Z_{1}^{n},\ldots,Z_{d}^{n}\right) \tol g\left(N_1,\ldots, N_d\right)
\]
where $N_1, \ldots, N_d$ are i.i.d.\ standard normal distributed random variables. 
\item[(b)] 
Let $\al \in (0,1)$ and let $k_\al$ denote the solution to the equation 
\[
\sum_{j=1}^q (-1)^{j+1} \binom{q}{j} \left(1-jk_\al \right)_+^{q-1} = \al.
\]
Then the test 
\begin{align} \label{test}
\vp(x) = \mathds{1}_{\left\{  f\left(Z_{1}^{n},\ldots,Z_{d}^{n}\right) > k_\al \right\}}
\end{align}
has asymptotic level $\al$ under $H_0$. 
\end{itemize}
\end{corollary}

When trying to understand the asymptotic behaviour of $f(Z_1^n, \ldots, Z_d^n)$ under the alternative the simplest case is when the limits $\left(e_{1},\ldots,e_{d}\right)$ are not in $A$. 

\begin{theorem} \label{thmaltA}
Let $(e_1, \ldots, e_d) \notin A$. Then 
\[
f\left(Z_{1}^{n},\ldots,Z_{d}^{n}\right) \pn g\left(e_1, \ldots, e_d\right).
\]
\end{theorem}

Theorem \ref{thmaltA} shows that the test from Corollary \ref{thmnull}(b) typically has power in the case where $(e_1, \ldots, e_d) \notin A$ holds, and the power increases the farther $g\left(e_1, \ldots, e_d\right)$ is away from the critical value $k_\al$. It hence remains to check in with situations $(e_1, \ldots, e_d) \notin A$ actually holds. We will start with a negative result regarding the case where $r$ and $d$ do not have common divisors. 

\begin{theorem} \label{thmalt}
Let $r$ and $d$ do not have common divisors. Then $(e_1, \ldots, e_d) \in A$ and
	\[
f\left(Z_{1}^{n},\ldots,Z_{d}^{n}\right) \tol g\left(N_1,\ldots, N_d\right)
\]
where $N_1, \ldots, N_d$ are i.i.d.\ standard normal distributed random variables. 
\end{theorem}

Theorem \ref{thmalt} indicates that we cannot expect the test to have any power when $r$ and $d$ do not have common divisors as in this case the asymptotic behaviour of $f\left(Z_{1}^{n},\ldots,Z_{d}^{n}\right)$ is the same as under the null hypothesis. 

For a full understanding of the asymptotics of the test, we need to discuss the situation where $r$ and $d$ \textit{do} have common divisors. The following theorem distinguishes the two cases of $r=2$ and $r \ge 3$.

\begin{theorem} \label{thmalt2}
Let $r$ and $d$ do have common divisors. 
\begin{itemize}
	\item[(a)] For $r=2$ we have $(e_1, \ldots, e_d) \in A$ always and 
	\[
f\left(Z_{1}^{n},\ldots,Z_{d}^{n}\right) \tol g\left(v_1 N_1, v_2 N_2, v_1 N_3, \ldots v_2 N_d\right)
\]
where $N_1, \ldots, N_d$ are i.i.d.\ standard normal distributed random variables.
	\item[(b)] Let $r \ge 3$ and let $b$ denote the largest common divisor of $r$ and $d$. If 
	\[
\sum_{k=1}^r e_{k}\exp\left(-\frac{2 \pi ik}{b}\right)\left(\left\lfloor \frac{d-k}{r}\right\rfloor +1\right)\neq0
\]
holds, then $(e_1, \ldots, e_d) \notin A$ and
\[
f\left(Z_{1}^{n},\ldots,Z_{d}^{n}\right) \pn g\left(e_1, \ldots, e_d\right).
\]
\end{itemize}
\end{theorem}

Theorem \ref{thmalt2}(a) gives a limiting distribution which, to the best of the authors' know\-ledge, has not yet been fully explored in the literature. Closest in spirit is a series of papers by Herbst in the 1960s in which, among other things, the distribution of the plain periodograms of normal distributions with a periodic variance is discussed. See e.g.\ Section 2 in \cite{herbst1963} where it is shown that these periodograms are not i.i.d.\ any more and become correlated in particular. This proves that the critical values obtained in Corollary \ref{thmnull}(b) most likely do not equal the ones for the limiting distribution in Theorem \ref{thmalt2}(a), although they may still be close in practical situations. 

Theorem \ref{thmalt2}(b), on the other hand, shows that the test often has significant power if $r \ge 3$ in combination with $r$ and $d$ having a common divisor. 

\begin{remark} \label{remunsys}
\mbox{} \vspace{-.2cm}
\begin{itemize}
	\item[(a)] While our model assumption is that $p_1, p_2, \ldots$ is periodic with an unknown frequency $r$, a natural question is to ask what happens in the case of a non-periodic choice of the sequence $p_1, p_2, \ldots$. While it is difficult to make a general statement, it is reasonable to assume a behaviour close to the one under the null hypothesis for most unsystematic choices of the success probabilities. 

For example, suppose that the sequence $p_1, p_2, \ldots$ is chosen i.i.d.\ according to some random variables $V_i$ with values in $[0,1]$, and with the conditional distribution of $Y_i$ given $V_i = p_i$ being independent $B(1,p_i)$ again. If the $(V_i,Y_i)$, $1 \le i \le n$, are then independent as before, it is obvious that 
\[
e_i = \lim_{n \to \infty} \E[Z_i^n] = \lim_{n \to \infty} \frac{1}{\left\lfloor \frac{n}{d}\right\rfloor }\sum_{k=0}^{\left\lfloor \frac{n}{d}\right\rfloor -1}\mathbb{E}\left[Y_{i+kd}\right] = \E[Y_1] = \mathbb{E}\left[\E\left[Y_{1}|V_1\right]\right] = \E\left[V_1\right]
\]
holds for any fixed $1 \le i \le d$. Hence $(e_1, \ldots, e_d) \in A$, and it is easy to deduce that the asymptotic behaviour of the corresponding test statistic is according to Corollary \ref{thmnull}(a) as well. 
	\item[(b)] A close inspection of the proofs of our results shows that most of them hold true for more general forms of categorical data than just binary ones. Under the null hypothesis of the $Y_i$ being i.i.d.\ with values in $\{0,1,\ldots,K-1\}$ the proofs go through with only minor modifications. What becomes more difficult, however, is to state general results about the behaviour under the alternative as our test statistic is tailored to identify periodic changes in the \textit{expectation} of the $Y_i$, and these expectations may obviously be constant even if the \textit{distribution} of the $Y_i$ changes periodically. For binary data, there is a one-to-one relation between a constant expectation and a constant $p_i$.
	\item[(c)] If we deviate from the assumption of independence of the $Y_i$ then the statement of Corollary \ref{thmnull} fails to hold. Even in the case of weak dependence, the random variables $Z_1^n, \ldots, Z_d^n$ will in general be dependent, and this dependence carries over to the limiting normal distribution $(N_1, \ldots, N_d)$. Similar to the context of Theorem \ref{thmalt2}(a), the behaviour of $g\left(N_1, \ldots, N_d\right)$ in this case is not known in the literature, even though in practice it might again not deviate much from the i.i.d.\ situation.
\end{itemize}

\end{remark}

\section{Finite sample properties} \label{sec:sim}
\def\theequation{3.\arabic{equation}}
\setcounter{equation}{0}

\subsection{Simulation study} \label{subsec:sim}

\begin{table}
\begin{center}
\begin{tabular}{|c|c|c|c|c|c|c|c|c|c|}
\hline 
$p_{1}$ & 0.1 & 0.2 & 0.3 & 0.4 & 0.5 & 0.6 & 0.7 & 0.8 & 0.9\tabularnewline
\hline 
 & 0.0487 & 0.0496 & 0.0504 & 0.0508 & 0.0498 & 0.0507 & 0.0498 & 0.0506 & 0.0497\tabularnewline
\hline 
\end{tabular}
\caption{Empirical level under the null hypothesis for various values of the true $p_1$.}
\label{tab1}
\end{center}
\end{table}

The following simulation study discusses the finite sample behaviour of the proposed test statistic in Corollary \ref{thmnull}(b). Unless stated otherwise, we work with $n=1200$ and $d=60$, and we use at least 20,000 repetitions to compute the empirical level and the empirical power, respectively. In all cases, the test is constructed using $k_\al$ for $\al = 0.05$. In practice, solving 
\[
\sum_{j=1}^q (-1)^{j+1} \binom{q}{j} \left(1-jx \right)_+^{q-1} = \al
\]
explicitly is often impossible, and we will follow the guideline given in \cite{fisher1929} and use the approximation of $k_\al$ via the last summand in the sum above only, i.e.\ we choose $k_\al$ in practice as the solution to
\begin{align} \label{kappr}
q\left(1-x \right)^{q-1} = \al.
\end{align}

Table \ref{tab1} proves that the level is kept in all situations, with empirical values being very close to $5\%$ always.

\begin{table}
\begin{center}
\begin{tabular}{|c|c|c|c|c|c|c|c|c|c|}
\hline 
$r$ & 2 & 3 & 4 & 5 & 6 & 7 & 8  & 9 & 10\tabularnewline
\hline 
 & 0.0518 & 0.0497 & 0.0530 & 0.0506 & 0.0560 & 0.0493 & 0.0495 & 0.0510 & 0.0684\tabularnewline
\hline 
11 & 12 & 13 & 14 & 15 & 16 & 17 & 18 & 19 & 20\tabularnewline
\hline 
 0.0497 & 0.0836 & 0.0484 & 0.0509 & 0.1350 & 0.0536 & 0.0503 & 0.0528 & 0.0507 & 0.2939 \tabularnewline
\hline 
 21 & 22 & 23 & 24 & 25 & 26 & 27 & 28 & 29 & 30\tabularnewline
\hline 
 0.0492 & 0.0518 & 0.0500 & 0.0903 & 0.0505 & 0.0504 & 0.0480 & 0.0498 & 0.0541 & 0.7473\tabularnewline
\hline 
\end{tabular}
\caption{Rejection level under the alternative of a periodic $p_{1},\ldots, p_{r}$ such that $p_{i+1}-p_{i}=0.01$ and $\frac{1}{r} \sum_{k=1}^r p_{i}=0.5$.}
\label{tab2}
\end{center}
\end{table}

For the alternative, we start with a specific model in which we discuss two different scenarios. For the first scenario, we choose $p_{1},\ldots, p_{r}$ in a periodic way and such that $p_{i+1}-p_{i}=0.01$ for all $1\leq i\leq r-1$, where the mean $\frac{1}{r} \sum_{k=1}^r p_{i}$ always equals $0.5$.

Because of the small differences in the probabilities, particularly for a small $r$, the alternatives in Table \ref{tab2} are difficult to detect. The situation improves for larger $r$, where in particular $r=20$ and $r=30$ (which have a lot of common divisors with $d=60$) give decent numbers. Note also that the power is close to $5\%$ in all cases where $r$ and $d$ do not have any common divisors, well in accordance with Theorem \ref{thmalt}. 

In Table \ref{tab3} the periodic behaviour is more pronounced, with the general setting being as in Table \ref{tab2} but now with $p_{i+1}-p_{i}=0.02$. Hence, we see an improvement across the board, and in particular $r=10, 12, 15, 24$ now have a large power as well. Additional simulations with $p_{i+1}-p_{i}=0.05$ (and then necessarily $r \le 20$) not reported here confirm these findings, and in that case essentially every alternative with $r \ge 3$ having a common divisor with $d$ is detected. 

\begin{table}
\begin{center}
\begin{tabular}{|c|c|c|c|c|c|c|c|c|c|}
\hline 
$r$ & 2 & 3 & 4 & 5 & 6 & 7 & 8  & 9 & 10\tabularnewline
\hline 
 & 0.0520 & 0.0550 & 0.0597 & 0.0779 & 0.0980 & 0.0482 & 0.0619 & 0.0551 & 0.3204 \tabularnewline
\hline 
11 & 12 & 13 & 14 & 15 & 16 & 17 & 18 & 19 & 20\tabularnewline
\hline 
 0.0488 & 0.4966 & 0.0518 & 0.0491 & 0.7550 & 0.0638 & 0.0503 & 0.1014 & 0.0504 & 0.9685\tabularnewline
\hline 
 21 & 22 & 23 & 24 & 25 & 26 & 27 & 28 & 29 & 30\tabularnewline
\hline 
0.0577 & 0.0521 & 0.0519 & 0.5326 & 0.0802 & 0.0507 & 0.0561 & 0.0654 & 0.1160 & 1\tabularnewline
\hline 
\end{tabular}
\caption{Rejection level under the alternative of a periodic $p_{1},\ldots, p_{r}$ such that $p_{i+1}-p_{i}=0.02$ and $\frac{1}{r} \sum_{k=1}^r p_{i}=0.5$.}
\label{tab3}
\end{center}
\end{table}

Table \ref{tab4}, on the other hand, discusses the situation where not the difference between two $p_i$ is kept constant in $r$ but in this case it is the smallest and the largest probability. Precisely, we set $p_1=0.4$ and $p_r=0.6$ and choose $p_{1},\ldots, p_{r}$ periodically. In comparison with the previous tables, it is now also the very small period lengths ($r \ge 3$ and such that $r$ divides $d$, so primarily $3 \le r \le 6$) which are often detected. This is clearly connected with the fact that the periodic behaviour for these large frequencies is now much more pronounced, particularly when compared with the model underlying Table \ref{tab2}. 

\begin{table}
\begin{center}
\begin{tabular}{|c|c|c|c|c|c|c|c|c|c|}
\hline 
$r$ & 2 & 3 & 4 & 5 & 6 & 7 & 8  & 9 & 10\tabularnewline
\hline 
 & 0.0492 & 0.9750 & 0.8362 & 0.6762 & 0.5988 & 0.0521 & 0.0983 & 0.0629 & 0.4306\tabularnewline
\hline 
11 & 12 & 13 & 14 & 15 & 16 & 17 & 18 & 19 & 20\tabularnewline
\hline 
0.0498 & 0.3838  & 0.0504  & 0.0519  & 0.3398  & 0.0533 & 0.0488 & 0.0575  & 0.0526  & 0.3411\tabularnewline
\hline 
 21 & 22 & 23 & 24 & 25 & 26 & 27 & 28 & 29 & 30\tabularnewline
\hline 
0.0489 & 0.0504 & 0.0495 & 0.0726 & 0.0499 & 0.0500 & 0.0529 & 0.0490 & 0.0505 & 0.3172 \tabularnewline
\hline 
\end{tabular}
\caption{Rejection level under the alternative of a periodic $p_{1},\ldots, p_{r}$ such that $p_1=0.4$ and $p_r=0.6$.}
\label{tab4}
\end{center}
\end{table}

Finally, Table \ref{tab5} shows the case where $p_{1},\ldots,p_{r}$ are given as the function values of $x\mapsto\frac{2}{5}\sin(x)+\frac{1}{2}$ with $p_{1}=\frac{2}{5}\sin(0)+\frac{1}{2}=\frac{1}{2}$, $p_{r}=\frac{2}{5}\sin(4\pi)+\frac{1}{2}=\frac{1}{2}$ and the arguments $x$ chosen equidistantly between $0$ and $4\pi$. That $r=3$ and $r=5$ are not detected is simply for symmetry reasons as in this case all success probabilities are actually the same, so we are in the case $r=1$ again.

\begin{table}
\begin{center}
\begin{tabular}{|c|c|c|c|c|c|c|c|c|c|}
\hline 
$r$ & 2 & 3 & 4 & 5 & 6 & 7 & 8 & 9 & 10\tabularnewline
\hline 
& 0.0507 & 0.0505 & 1 & 0.0496 & 1 & 0.0513 & 1 & 1 & 1\tabularnewline
\hline 
\end{tabular}
\caption{Rejection level under the alternative of a periodic $p_{1},\ldots, p_{r}$ chosen according to the function $f$.}
\label{tab5}
\end{center}
\end{table}

We end this simulation study with a simple example that neither fits into $H_{0}$ nor into $H_{1}$, namely non-periodic success probabilities. In this case we have set $n=120$ and $d=12$, and we have chosen $p_{i}=a_i^{-1}$ where $a_i$ denotes the $i$th digit in the decimal representation of $\pi$, i.e.\ $p_{1},p_{2},p_{3},\ldots,p_{120}=0.1,0.4,0.1,\ldots,0.8$. In this case our simulations give $\mathbb{P}(f(Z_1^n, \ldots, Z_d^n)\geq k_{0.05})=0.0529$, a behaviour which is very much in line with the content of Remark \ref{remunsys}(a).

\subsection{Data analysis}
We have applied our test statistic to two series of binary data connected with floods at the rivers Elbe and Oder. Both series have already been discussed in \cite{mudetal2003} over a time period from the 11th to the end of the 20th century, and it was argued in the paper that no upward trend in the occurrence of floods can be detected. This makes it reasonable to ask whether a periodic behaviour exists instead. A natural candidate for a seasonality is obviously given by $r=12$, but as there is a tendency to both summer and winter floods $r=6$ or even other seasonalities might make sense. To ensure for a somewhat constant data quality over time, we have decided to work with floods from 1903 to 2002 only, as there are less and less floods recorded the further we go back in time. We are hence working with a set of data of length $n=1200$ and we can use the same setup as in Section \ref{subsec:sim} and use $d=60$ and, hence, $q=29$. 

For the Elbe, the series shows a list of 40 floods over the period of 1200 months. In cases where the duration of a flood was longer than a month, we have decided to incorporate the flood only once and to use the first of the respective months as the data point, while treating all other months during that period as showing no flood. Since the relative frequency of floods is low we have first conducted an additional simulation study and checked the behaviour of the test for a constant success probability of $p_1=0.03$, a number close to the average number of floods per month. In the same setup as in Table \ref{tab1} we obtained an empirical level of 0.0442 which is slightly more conservative than for the other cases in Table \ref{tab1} but still very reasonable. For the Oder we have 79 floods in the corresponding time period and proceeded in the same way, refraining from additional simulations as the relative amount of floods is now reasonably close to the already tested $p_1=0.1$. 

For $q=29$ and with the same approximation as in \eqref{kappr} we obtain a critical value of $k_{0.05}=0.2033$. For the Elbe data the value of the test statistic equals 0.1356 (with the maximising periodogram at $j=27$) which corresponds to a $p$-value of 0.4900 and suggests that no periodic behaviour exists. The Oder data behaves similarly overall. Here, the value of the test statistic is 0.1414 with a p-value of 0.4062 and a maximising periodogram at $j=5$. Both findings come with a grain of salt, however, as we know in particular from Tables \ref{tab2} and \ref{tab3} that the power of the test is not that high for the candidate periodicities $r=6$ and $r=12$ when the difference in success probabilities is relatively small. Given that the number of floods is relatively low to begin with, it is therefore a priori not surprising that the test indeed opts to accept the hypothesis of a constant probability. 

\section{Proofs} \label{sec:proof}
\def\theequation{4.\arabic{equation}}
\setcounter{equation}{0}

\subsection{Proof of Theorem \ref{thmei}}
As a first step, we prove 
\[
Z_{i}^{n}-\frac{1}{\left\lfloor \frac{n}{d}\right\rfloor }\sum_{k=0}^{\left\lfloor \frac{n}{d}\right\rfloor -1}\mathbb{E}\left[Y_{i+kd}\right] = o_\P(1)
\]
which by definition of $Z_i^n$ follows from
\begin{align*}
&\mathbb{E}\left[\left(\frac{1}{\left\lfloor \frac{n}{d}\right\rfloor }\sum_{k=0}^{\left\lfloor \frac{n}{d}\right\rfloor -1}\left(Y_{i+kd}-\mathbb{E}\left[Y_{i+kd}\right]\right)\right)^{2}\right] 
\\ =& \frac{1}{\left(\left\lfloor \frac{n}{d}\right\rfloor\right)^2}\sum_{k=0}^{\left\lfloor \frac{n}{d}\right\rfloor -1}\mathrm{Var}\left(Y_{i+kd}\right)
= \frac{1}{\left(\left\lfloor \frac{n}{d}\right\rfloor\right)^2}\sum_{k=0}^{\left\lfloor \frac{n}{d}\right\rfloor -1}p_{i+kd}(1-p_{i+kd}) \le \frac{1}{4\left\lfloor \frac{n}{d}\right\rfloor} \longrightarrow 0
\end{align*}
where we have used the independence of the $Y_i$ as well as the fact $p_{i+kd}(1-p_{i+kd})$ is bounded by $1/4$. 

Hence, it remains to show
\[
\underset{n\rightarrow\infty}{\lim}\frac{1}{\left\lfloor \frac{n}{d}\right\rfloor }\sum_{k=0}^{\left\lfloor \frac{n}{d}\right\rfloor -1}\mathbb{E}\left[Y_{i+kd}\right]=e_{i}
\]
for which we first rewrite 
\begin{align*}
&\sum_{k=0}^{\lfloor \frac{n}{d}\rfloor -1}\mathbb{E}\left[Y_{i+kd}\right]=\sum_{k=0}^{\lfloor \frac{n}{d}\rfloor -1}p_{i+kd} =\sum_{j=0}^{r-1}\sum_{m=0}^{\lfloor \frac{\lfloor \frac nd \rfloor}{r}  \rfloor -1}p_{i+(j+mr)d} 
+ \sum_{s=r \lfloor \frac{\lfloor \frac nd \rfloor}{r}  \rfloor}^{\lfloor \frac{n}{d}\rfloor -1} p_{i+sd}.
\end{align*}
Due to the trivial bound $r \lfloor \frac{\lfloor \frac nd \rfloor}{r}  \rfloor \ge r (\frac{\lfloor \frac nd \rfloor}{r} -1)$, the second sum contains at most $r$ terms. Hence, 
\begin{align*}
\frac{1}{\left\lfloor \frac{n}{d}\right\rfloor } \sum_{s=r \lfloor \frac{\lfloor \frac nd \rfloor}{r}  \rfloor}^{\lfloor \frac{n}{d}\rfloor -1} p_{i+sd} \longrightarrow 0.
\end{align*}
For the double sum we can obviously write 
\begin{align*}
 \sum_{j=0}^{r-1}\sum_{m=0}^{\lfloor \frac{\lfloor \frac nd \rfloor}{r}  \rfloor -1}p_{i+(j+mr)d} = \frac 1r \sum_{j=0}^{r-1} r \sum_{m=0}^{\lfloor \frac{\lfloor \frac nd \rfloor}{r}  \rfloor -1}p_{i+(j+mr)d} = \frac 1r \sum_{j=0}^{r-1} r \sum_{m=0}^{\lfloor \frac{\lfloor \frac nd \rfloor}{r}  \rfloor -1}p_{i+jd}
\end{align*}
by periodicity of the $p_i$. The terms in the latter sum are in fact independent of $m$, and setting $\langle x \rangle = x - \lfloor x \rfloor$ we obtain
\begin{align*}
 &\frac{1}{\left\lfloor \frac{n}{d}\right\rfloor }\sum_{j=0}^{r-1}\sum_{m=0}^{\lfloor \frac{\lfloor \frac nd \rfloor}{r}  \rfloor -1}p_{i+(j+mr)d} = \frac 1r \sum_{j=0}^{r-1} \frac{r}{\left\lfloor \frac{n}{d}\right\rfloor } \left \lfloor \frac{\lfloor \frac nd \rfloor}{r}  \right \rfloor p_{i+jd} \\ =& \frac 1r \sum_{j=0}^{r-1} \frac{r}{\left\lfloor \frac{n}{d}\right\rfloor } \left( \frac{\lfloor \frac nd \rfloor}{r}  - \left \langle \frac{\lfloor \frac nd \rfloor}{r} \right \rangle \right) p_{i+jd} = e_i - \frac 1r \sum_{j=0}^{r-1} \frac{r}{\left\lfloor \frac{n}{d}\right\rfloor }  \left \langle \frac{\lfloor \frac nd \rfloor}{r} \right \rangle p_{i+jd}.
\end{align*}
The claim then follows from 
\begin{align*}
\frac 1r \sum_{j=0}^{r-1} \frac{r}{\left\lfloor \frac{n}{d}\right\rfloor }  \left \langle \frac{\lfloor \frac nd \rfloor}{r} \right \rangle p_{i+jd} \le \frac{1}{\left\lfloor \frac{n}{d}\right\rfloor }  \sum_{j=0}^{r-1} p_{i+jd} \longrightarrow 0
\end{align*}
using $\langle x \rangle \le 1$. \qed

\subsection{Proof of Lemma \ref{lemei}}

To prove part (a) let $j\in\mathbb{N}$ and $1\leq k\leq r$. Then we have
\[
p_{i+jr+kd}=p_{i+kd}
\]
by $r$-periodicity of the $p_{i}$, so we obtain
\[
e_{i+jr}=\frac{1}{r}\sum_{k=0}^{r-1}p_{i+jr+kd}=\frac{1}{r}\sum_{k=0}^{r-1}p_{i+kd}=e_{i}.
\]
%Part (b) follows in a similar way. For $m,j\in\mathbb{N}$ and $0\leq s\leq r-1$ it is clear that
%\[
%p_{i+(s+jr)d}=p_{i+sd+jrd}=p_{i+sd}.
%\]
%It is then easy to conclude 
%\[
%\sum_{k=0}^{mr-1}p_{i+kd}=\sum_{j=0}^{m-1}\sum_{s=0}^{r-1}p_{i+(s+jr)d}=\sum_{j=0}^{m-1} \sum_{s=0}^{r-1}p_{i+sd}=m\sum_{s=0}^{r-1}p_{i+sd}
%\]
%and hence
%\[
%\frac{1}{mr}\sum_{s=1}^{mr-1}p_{i+sd} =\frac{1}{mr}m\sum_{s=0}^{r-1}p_{i+sd}=\frac{1}{r}\sum_{s=0}^{r-1}p_{i+sd}=e_{i}.
%\]

For part (b), our strategy of proof is to show that, regardless of the choice of $i$, all $r$ representatives of the residue classes (i.e.\ any of the $p_{1},\ldots,p_{r}$) appear once within the sum $\sum_{k=0}^{r-1}p_{i+kd}$, from which it follows that the $e_i$ are constant in $i$. We will give a proof by contraction which relies on the simple fact that $r$ and $d$ having no common divisors is equivalent to the least common multiple of $r$ and $d$ being $rd$. 

To derive the contradiction let us assume that there exist $1\leq k < \ell\leq r$ such that $i+kd=i+\ell d\ (\mathrm{mod}\ r)$. Then there exists $j\in\mathbb{N}$ such that $i+kd+jr=i+\ell d$, and this implies first $kd+jr=\ell d$ and then $jr=(\ell-k)d$. Because of $\ell-k<r$ this is obviously a least common multiple smaller than $rd$, hence there is a contradiction.

Finally, for part (c) let us assume the existence of a common divisor $m\in\mathbb{N}$ with $1<m\le r$, i.e.\ let us assume $r=m\ell$ and $d=mk$ for integers $k,\ell\in\mathbb{N}$. In particular, $\ell<r$, and we conclude that
\[
\ell d=\frac{r}{m} m k=rk.
\]
Hence
\[
p_{i+\ell d}=p_{i+rk}=p_{i}
\]
for any $0 \le i \le r-1$, and hence each sum $\sum_{k=0}^{r-1}p_{i+kd}$ only contains summands of the form $p_{i},p_{i+d},\ldots,p_{i+(\ell-1)d}$. Since $\ell<r$
this proves that not all representatives of the $r$ residue classes are contained in $e_i$, and hence the $e_{i}$ are in general not constant. \qed

\subsection{Proof of Lemma \ref{lemeiA}}
Let first $\left(Z_{1}^n,\ldots, Z_{d}^n\right) \notin A$. Then a straight forward calculation gives
\begin{align*}
&g\left(c_n(Z_{1}^{n}-e_{1}), \ldots,c_n(Z_{d}^{n}-e_{d})\right)\\
%=&\max_{j=1, \ldots, q} \frac{I_{c_n\cdot(Z-e)}^n(\om_j)}{\sum_{m=1}^q I_{c_n\cdot(Z-e)}^n(\om_m)}   
=&\frac{\max_{j=1, \ldots, q}  \frac 1d \left| \sum_{\ell=1}^d c_n(Z_{\ell}^{n}-e_{\ell}) \exp(-i \ell \om_j) \right |^2}{\sum_{m=1}^d \frac 1d \left| \sum_{\ell=1}^d c_n(Z_{\ell}^{n}-e_{\ell}) \exp(-i \ell \om_m) \right |^2}
\\=& \frac{\max_{j=1, \ldots, q}  \frac 1d \left| \sum_{\ell=1}^d (Z_{\ell}^{n}-e_{\ell}) \exp(-i \ell \om_j) \right |^2}{\sum_{m=1}^d \frac 1d \left| \sum_{\ell=1}^d (Z_{\ell}^{n}-e_{\ell}) \exp(-i \ell \om_m) \right |^2}
\\=& \frac{\max_{j=1, \ldots, q} \frac 1d \left| \sum_{\ell=1}^d Z_{\ell}^{n} \exp(-i \ell \om_j) \right |^2}{\sum_{m=1}^d \frac 1d \left| \sum_{\ell=1}^d Z_{\ell}^{n} \exp(-i \ell \om_m) \right |^2} = g\left(Z_{1}^{n},\ldots,Z_{d}^{n}\right)
\end{align*}
because of 
\[
\sum_{\ell=1}^d e_{\ell} \exp(-i \ell \om_j) =0
\] 
for each $1\leq j\leq q$ by assumption. 

By definition of $f$ it hence remains to prove 
\[
\left(Z_{1}^{n},\ldots,Z_{d}^{n}\right) \in A \: \Longleftrightarrow \: \left(c_n(Z_{1}^{n}-e_{1}), \ldots, c_n(Z_{d}^{n}-e_{d})\right) \in A
\]
which again follows easily from the assumption and the definition of $A$.
\qed

\subsection{Proof of Theorem \ref{thmgen}}
We will use the version of the continuous mapping theorem given as Theorem 1.3.6 in \cite{vanwel1996}. To this end, note first that
\[
f(x_1, \ldots, x_d) = g(x_1, \ldots, x_d) \mathds{1}_{\left\{ \left(x_1, \ldots, x_d\right)\notin A\right\}}
\]
is a function which is continuous outside of the set $A$, and that Lemma \ref{lemeiA} gives
\[
f\left(Z_{1}^{n},\ldots,Z_{d}^{n}\right) =  f\left(c_n(Z_{1}^{n}-e_{1}), \ldots, c_n(Z_{d}^{n}-e_{d})\right).
\]
Setting 
\[
c_n = \sqrt{\left\lfloor \frac{n}{d}\right\rfloor}
\]
the central limit theorem together with Slutsky's theorem and Remark \ref{remei} gives 
\[
\left(c_n(Z_{1}^{n}-e_{1}), \ldots, c_n(Z_{d}^{n}-e_{d})\right) \tol \left(v_1 N_1, \ldots, v_d N_d \right).
\]
Thus, it remains to show 
\begin{align} \label{nullA}
\P\left(\left(v_1 N_1, \ldots, v_d N_d \right) \in A \right) =0
\end{align}
in order to deduce 
\[
f\left(Z_{1}^{n},\ldots,Z_{d}^{n}\right) \tol f\left(v_1 N_1, \ldots, v_d N_d\right) \stackrel{d}{=} g\left(v_1 N_1, \ldots, v_d N_d\right).
\]
To prove (\ref{nullA}) note that e.g.\ Theorem I.9.2 and Example III.3(f) in \cite{feller1971} prove, among other things, that $g\left(N_1, \ldots, N_d\right)$ is well defined almost surely. By the Radon-Nikod\'ym theorem it is clear that $\left(v_1 N_1, \ldots, v_d N_d \right)$ and $\left(N_1, \ldots, N_d \right)$ share the same null sets, which gives (\ref{nullA}). 

Finally, by definition $g\left(x_1, \ldots, x_d\right) = g\left(yx_1, \ldots, yx_d\right)$ holds on the set $A$ for any choice of $y>0$. \qed

\subsection{Proof of Corollary \ref{thmnull}}
For part (a) one uses the fact that $e_i=p_1$ and $v_i = p_1(1-p_1)$, $1 \le i \le d$, are constant which follows directly from Remark \ref{remei}. Then Theorem \ref{thmgen} immediately gives the claim since $\left(e_{1},\ldots,e_{d}\right) \in A$ by the geometric sum argument 
\[
\sum_{\ell=1}^d \exp(-i \ell \om_j)= \sum_{\ell=1}^d \exp\left(-\frac{2 \pi i \ell j}d\right) = \sum_{\ell=0}^{d-1} \exp\left(-\frac{2 \pi i \ell j}d\right) = \frac{1-\exp\left(-{2 \pi i j}\right)}{1-\exp\left(-\frac{2 \pi i j}d\right)} = 0.
\]
Note that the denominator does not become zero because of $1 \le j \le q$. Critical values for the corresponding level $\al$ test in (b) can then be derived using (\ref{limxi}). \qed

\subsection{Proof of Theorem \ref{thmaltA}}
Here we use the extended continuous mapping theorem as stated in Theorem 1.11.1 in \cite{vanwel1996}. In their notation we set $\D_n = \R^d$ and $\D_0 = \R^d \backslash A$. Since $\D_0$ is an open set it is clear that each sequence $(x_n)$ in $\D$ with $x_n \longrightarrow x \in \D_0$ satisfies $x_n \in \D_0$ eventually. In particular, for any such sequence we have $f(x_n) \longrightarrow g(x)$  by definition of $f$ and $g$. Since 
\[ 
\left(Z_{1}^{n},\ldots,Z_{d}^{n}\right) \pn \left(e_1, \ldots, e_d\right)
\]
holds as a consequence of Theorem \ref{thmei} and $\left(e_1, \ldots, e_d\right) \in \D_0$ by assumption we may conclude
\[
f\left(Z_{1}^{n},\ldots,Z_{d}^{n}\right) \pn g\left(e_1, \ldots, e_d\right)
\]
from the afore mentioned Theorem 1.11.1 in \cite{vanwel1996}. \qed

\subsection{Proof of Theorem \ref{thmalt}}
Note first that $\left(e_1, \ldots, e_d\right) \in A$ which, as in the proof of Corollary \ref{thmnull}, is an easy consequence of Lemma \ref{lemei}(b) and the definition of $A$. Theorem \ref{thmgen} then can be applied, and Remark \ref{remei} again gives 
\[
v_i = \underset{n\rightarrow\infty}{\lim}\Var\left(\sqrt{\left\lfloor \frac{n}{d}\right\rfloor} Z_i^n\right) = \frac{1}{r} \sum_{j=0}^{r-1} p_{i+jd}(1-p_{i+jd})
\]
for any $1 \le i \le d$. The argument from the proof of Lemma \ref{lemei}(b) then shows that the limit is in fact independent of the choice of $i$ and equals $\frac{1}{r} \sum_{k=0}^{r-1} p_{k}(1-p_{k})$. Theorem \ref{thmgen} concludes. \qed

\subsection{Proof of Theorem \ref{thmalt2}}
For part (a) note first that $d$ has to be even as it otherwise cannot have a common divisor with $r=2$. Also, we know from Lemma \ref{lemei}(a) that $(e_1, \ldots, e_d)$ is $r$-periodic. Now, let $j\in\{1,\ldots,q\}$ with $q=\frac d2 -1$. Then
\begin{align*}
\sum_{\ell=1}^d e_{\ell}\exp\left(-i\ell\omega_{j}\right)&=\sum_{k=1}^r\sum_{m=0}^{\frac dr -1}e_{k+mr}\exp\left(-i(k+mr)\omega_{j}\right)\\
%&=\sum_{k=1}^r\sum_{m=0}^{\frac dr -1}e_{k}\exp\left(-i(k+mr)\omega_{j}\right)\\ 
&=\sum_{k=1}^r\sum_{m=0}^{\frac dr -1}e_{k}\exp\left(-i(k+mr)\omega_{j}\right) \\
&=\sum_{k=1}^re_{k}\exp\left(-ik\omega_{j}\right)\sum_{m=0}^{\frac dr -1}\exp\left(-imr\omega_{j}\right).
%\\ &=\sum_{k=1}^r e_{k}\exp\left(-ik\omega_{j}\right)\sum_{m=1}^{\frac dr -1}\exp\left(-ir\omega_{j}\right)^{m}
\end{align*}
A geometric sum argument now gives 
\[
\sum_{m=0}^{\frac dr -1}\exp\left(-imr\omega_{j}\right) = \sum_{m=0}^{\frac dr -1}\exp\left(-\frac{2 \pi i mrj}d\right) = \frac{1-\exp\left(-{2 \pi i  m j}\right)}{1-\exp\left(-\frac{2 \pi i r j}d\right)} = 0.
\]
Again the denominator cannot become zero as in this case $\frac{rj}{d}\in\mathbb{Z}$ would have to hold. This is excluded in light of $1 \le j \le \frac d2-1$ and $r=2$. Theorem \ref{thmgen} now yields 
\[
f\left(Z_{1}^{n},\ldots,Z_{d}^{n}\right) \tol g\left(v_1 N_1,\ldots, v_d N_d\right),
\]
and by definition it is clear that each $v_i$ either equals $v_1$ or $v_2$, depending on whether $i$ is uneven or not. 

For part (b) we need to show $(e_1, \ldots, e_d) \notin A$ and then apply Theorem \ref{thmaltA}. We start with a simple auxiliary lemma and show that $\frac db \in \N$ gives a fine Fourier frequency, i.e.\ that 
$\frac{d}{b}\leq\left\lfloor \frac{d-1}{2}\right\rfloor$ holds. Clearly, $r \ge 3$ together with $r$ and $d$ having a common divisor implies both $b\ge 3$ and $d \ge 3$, and in particular it is enough to prove $\frac{d}{3}\leq\left\lfloor \frac{d-1}{2}\right\rfloor$.

In the case of $d$ being uneven, it is clear that $\left\lfloor \frac{d-1}{2}\right\rfloor =\frac{d-1}{2}$, and it is simple to deduce that $d \ge 3$ implies $\frac{d}{3}\leq\frac{d-1}{2}$. For an even $d$ we necessarily have $d \ge 6$, and another easy calculation shows that this implies   
\[\frac{d}{3}\leq\frac{d}{2} -1=\left\lfloor \frac{d-1}{2}\right\rfloor.
\]

Now, let $d$ be a multiple of $r$. In this case, the very same calculation as for part (a) gives 
\begin{align*}
\sum_{\ell=1}^d e_{\ell}\exp\left(-i\ell\omega_{\frac db}\right) =\sum_{k=1}^r e_{k}\exp\left(-ik\omega_{\frac db}\right)\sum_{m=0}^{\frac dr -1}\exp\left(-imr\omega_{\frac db}\right), 
\end{align*}
but now 
\[
\exp\left(-imr\omega_{\frac db}\right) = \exp\left(-\frac{2 \pi imr}b \right) = 1
\]
because $b$ divides $r$. Hence, 
\begin{align*}
\sum_{\ell=1}^d e_{\ell}\exp\left(-i\ell\omega_{\frac db}\right) &= \sum_{k=1}^r e_{k}\exp\left(-ik\omega_{\frac db}\right) \frac dr = \sum_{k=1}^r e_{k}\exp\left(-ik\omega_{\frac db}\right) \left(\left\lfloor \frac{d-k}{r}\right\rfloor +1\right),
\end{align*}
from which the claim can directly be derived. If $d$ is not a multiple of $r$, after 
\begin{align*}
\sum_{\ell=1}^d e_{\ell}\exp\left(-i\ell\omega_{\frac db}\right) =\sum_{k=1}^r e_{k}\exp\left(-ik\omega_{\frac db}\right)\sum_{m=0}^{\left\lfloor \frac{d-k}{r}\right\rfloor}\exp\left(-imr\omega_{\frac db}\right), 
\end{align*}
the same calculation as above works. \qed

\bibliographystyle{chicago}
\bibliography{literatur}
\end{document}